\documentclass [a4paper,11pt]{article}

\usepackage{amsthm}
\usepackage{amsmath,amssymb,latexsym,amsfonts,mathrsfs}
\usepackage[dvipdfmx]{graphicx}

\usepackage{color}

\usepackage{mathtools}
\usepackage{fancyhdr}
\usepackage[top=30truemm, bottom=30truemm, left=25truemm, right=25truemm]{geometry}

\newcommand{\ep}{\varepsilon}
\newcommand{\nn}{\nonumber}

\newcommand{\SCR}[1]{{\mathscr #1}}

\newcommand{\CAL}[1]{{\cal #1}
}

\newcommand{\J}[1]{\left\langle #1 \right\rangle}
\newcommand{\D}[1]{{\mathscr D}( #1 )}

\theoremstyle{definition}
\newtheorem{Thm}{{\bf Theorem}}[section]

\newtheorem{Lem}[Thm]{{\bf Lemma}}
\newtheorem{Prop}[Thm]{{\bf Proposition}}

\newtheorem{Ass}[Thm]{{\bf Assumption}}

\newtheorem{Rem}[Thm]{{\bf Remark}}

\newcounter{Exami}

\newcommand{\Proof}[2][Proof]{
\begin{proof}[{\bf #1}]
#2
\end{proof}
}

\begin{document}

%%%%%%%%%%%%%%%%%%%%%%%%%%%%%%%%%%%%%%%%%%%%%%%%%%%%%%%%%%%%%%%%%%%%%%%%%%%%%%%%%%%%%%%%%%%%%%%%%%%%%%%%%%%%%%%%%%%%%%%%%%%
%author
%%%%%%%%%%%%%%%%%%%%%%%%%%%%%%%%%%%%%%%%%%%%%%%%%%%%%%%%%%%%%%%%%%%%%%%%%%%%%%%%%%%%%%%%%%%%%%%%%%%%%%%%%%%%%%%%%%%%%%%%%%%
\begin{flushleft}
{\bf \Large Asymptotic behavior for nonlinear Schr\"{o}dinger equations with critical time-decaying harmonic potential.
} \\ \vspace{0.3cm}
by {\bf \large Masaki Kawamoto } \\ 
Department of Engineering for Production, Graduate School of Science and Engineering, Ehime University, 3 Bunkyo-cho Matsuyama, Ehime, 790-8577. Japan \\
Email: {kawamoto.masaki.zs@ehime-u.ac.jp} 
\end{flushleft}

%%%%%%%%%%%%%%%%%%%%%%%%%%%%%%%%%%%%%%%%%%%%%%%%%%%%%%%%%%%%%%%%%%%%%%%%%%%%%%%%%%%%%%%%%%%%%%%%%%%%%%%%%%%%%%%%%%%%%%%%%%%
%abst
%%%%%%%%%%%%%%%%%%%%%%%%%%%%%%%%%%%%%%%%%%%%%%%%%%%%%%%%%%%%%%%%%%%%%%%%%%%%%%%%%%%%%%%%%%%%%%%%%%%%%%%%%%%%%%%%%%%%%%%%%%%
\begin{center}
\begin{minipage}[c]{400pt}
{\bf Abstract}. {\small
Time-decaying harmonic oscillators yield dispersive estimates with weak decay, and change the threshold power of the nonlinearity between the short and the long range. In the non-critical case for the time-decaying harmonic oscillator, this threshold can be characterized by polynomial nonlinearities. However, in the critical case, it is difficult to characterize the threshold using only polynomial terms, and thus we use logarithmic nonlinear terms.  
}
\end{minipage}
\end{center}

\begin{flushleft}
{\bf Keywords}; Nonlinear scattering theory; Long-range scattering; Time-dependent harmonic oscillators; Time-dependent magnetic fields
\end{flushleft}
\begin{flushleft}
{\bf MSC}; primary 35Q55, secondly 35J10
\end{flushleft}
%%%%%%%%%%%%%%%%%%%%%%%%%%%%%%%%%%%%%%%%%%%%%%%%%%%%%%%%%%%%%%%%%%%%%%%%%%%%%%%%%%%%%%%%%%%%%%%%%%%%%%%%%%%%%%%%%%%%%%%%%%%
%Intro
%%%%%%%%%%%%%%%%%%%%%%%%%%%%%%%%%%%%%%%%%%%%%%%%%%%%%%%%%%%%%%%%%%%%%%%%%%%%%%%%%%%%%%%%%%%%%%%%%%%%%%%%%%%%%%%%%%%%%%%%%%%
\section{Introduction}
In this study, we consider nonlinear Schr\"{o}dinger (NLS) equations with time-dependent harmonic potentials:
\begin{align}\label{4}
\begin{cases}
i \partial _t u(t,x) - \left( - \Delta /2 + \sigma (t) |x|^2/2 \right) u(t,x) =  F(u(t,x)) u(t,x) , \\
u(0,x) = u_0 (x),
\end{cases}
\end{align}
where $(t,x) \in {\bf R} \times {\bf R}^n$, $n \in \{ 1,2,3\}$, and $F:{{\bf C} \to {\bf R}}$ is a nonlinear term to be defined later. Let
\begin{align*}
H_0 (t) := - \Delta /2 + \sigma (t) |x|^2/2. 
\end{align*} 
$ \sigma (t) |x|^2/2$ is called time-dependent harmonic potential with  coefficient $\sigma (t)$. We make the following assumption on $\sigma (t) \in L^{\infty} ({\bf R}_t)$. 
\begin{Ass} \label{A1}
 Let $\zeta _1 (t)$ and $\zeta _2 (t)$ be the fundamental solutions to the following equations: 
\begin{align} \label{1}
\zeta _j ''(t) + \sigma (t) \zeta _j (t) =0, \quad
\begin{cases}
\zeta _1 (0) = 1, \\
\zeta _1 ' (0) =0,
\end{cases}
\quad
\begin{cases}
\zeta _2 (0) = 0, \\
\zeta _2 ' (0) =1.
\end{cases}
\end{align} 
Then, there exist $c>0$, $r_0 \gg 1$, and $c_{1, \pm} ,c_{2, \pm} \notin \{ 0, \infty, - \infty \}$ such that for all $|t| > r_0$, we have
\begin{align*}
|\zeta _2 (t)| \geq c 
\end{align*} 
and 
\begin{align}\label{ad3}
\lim_{t \to \pm \infty} \frac{\zeta _1 (t)}{|t|^{1/2}} = c_{1, \pm}, \quad 
\lim_{t \to \pm \infty} \frac{\zeta _2 (t)}{ |t| ^{1/2}   \log |t| } = c_{2, \pm}.
\end{align}
Moreover, $\zeta_1 (t)$, $\zeta_2 (t)$, $\zeta_1' (t)$, and $\zeta_2' (t)$ are continuous functions in $t$.
\end{Ass}
\begin{Rem}
For simplicity, let us redefine $\sigma$ as
\begin{align} \label{2}
\sigma (t) = 
\begin{cases}
\sigma_{s} (t), & |t| \leq r_0, \\ 
1/(4t^2) , & |t| > r_0 
\end{cases} 
\end{align}
where $\sigma _s (t)$ is a bounded smooth function. Then, for $|t| \geq r_0$, the functions $y_1 (t) = |t|^{1/2} $ and $y_2 (t) = |t|^{1/2} \log |t|$ are the linearly independent solutions to $y'' (t) + \sigma (t) y(t) =0$. Hence, in \eqref{2}, by writing $\zeta _1 (t) = c_{1,1, \pm} y_1 (t) + c_{1,2 , \pm} y_2 (t) $ and $\zeta _2 (t) = c_{2,1, \pm } y_1 (t) + c_{2,2, \pm} y_2 (t) $, Assumption \ref{A1} implies $c_{1,1, \pm} = c_{1,\pm}$, $c_{1,2, \pm} = 0$, $c_{2,1, \pm} \in {\bf R}$, and $c_{2,2, \pm}=c_{2,\pm}$. We provide examples of $\sigma (t)$ satisfying Assumption \ref{A1} in \S{4}.   
\end{Rem}

We now let $U_0(t,s)$ be a propagator for $H_0 (t)$, that is, a family of unitary operators $\{ U_0(t,s)\}_{(t,s) \in {\bf R}^2} $ on $L^2({\bf R}^n)$ such that for all $t,s,\tau \in {\bf R}$, we have
\begin{align*}
& i \partial _t U_0(t,s) = H_0(t) U_0(t,s), \quad i \partial _s U_0(t,s) = -U_0(t,s) H_0(s), \\
& U_0(t, \tau) U_0(\tau, s) = U_0(t,s), \quad U_0(s,s) = \mathrm{Id}_{L^2({\bf R}^n)}, \quad U_0(t,s) \D{H_0(s)} \subset \D{H_0(s) }
\end{align*}
on $\D{H_0(s)}$. By using $\zeta _1 (t)$ and $\zeta _2 (t)$, the following \CAL{MDFM decomposition} can be obtained, see Korotyaev \cite{Ko} (also Carles \cite{Ca}, Kawamoto \cite{Ka2}, and Kawamoto-Muramatsu \cite{KM}).

\begin{Lem}
For $\phi \in \SCR{S}({\bf R}^n)$, let 
\begin{align*}
\left( \CAL{M}(\tau) \phi \right) (x) = e^{ix^2/(2 \tau)} \phi (x), \quad
\left(
\CAL{D}(\tau) \phi
\right) (x) = \frac{1}{(i \tau)^{n/2}} \phi (x/ \tau).
\end{align*}
Then, the following MDFM decomposition holds:
\begin{align*}
U_0(t,0) = \CAL{M} \left(  \frac{\zeta _2(t)}{\zeta _2 '(t)} \right) \CAL{D} (\zeta _2 (t)) \SCR{F} \CAL{M} \left( \frac{\zeta _2 (t)}{ \zeta _1 (t)} \right)
\end{align*}
\end{Lem}
%Hereafter we use the notation 
%\begin{align*}
%\CAL{M}_1 (t) = \CAL{M} \left(  \frac{\zeta _2(t)}{\zeta _2 '(t)} \right), \quad \CAL{M}_2 (t) =\CAL{M} \left( \frac{\zeta _2 (t)}{ \zeta _1 (t)} \right).
%\end{align*}
By the decomposition formula in \cite{Ko}, we have the following dispersive estimates: \\ 
For $\phi \in L^1 ({\bf R}^n)$, 
\begin{align*}
\left\| 
U_0(t,0) \phi 
\right\| _{\infty} \leq C |\zeta _2 (t) | ^{-n/2} \left\| \phi \right\|_1
\end{align*}
and 
\begin{align*}
\left\| 
U_0(t,s) \phi 
\right\| _{\infty} \leq C |\zeta _1 (s) \zeta _2 (t) - \zeta _1 (t) \zeta _2 (s) | ^{-n/2} \left\| \phi \right\|_1,
\end{align*}
where $\| \cdot  \|_{q} $, $1 \leq q \leq \infty$, denotes the norm on $L^q ({\bf R}^n)$. Hence, we immediately obtain the following proposition. 
\begin{Prop} Under Assumption \ref{A1}, for all $|t| > r_0$, there exists $C_{} >0$ such that 
\begin{align}\label{3}
\left\| 
U_0(t,0) \phi 
\right\| _{\infty} \leq C | t | ^{-n/4} (\log |t|)^{-n/2} \left\| \phi \right\|_1. 
\end{align}
\end{Prop}
We now introduce the definition of nonlinearity $F$. Consider the case where $\sigma (t) \equiv 0$ and $F(u(t,x)) = |u(t,x)|^{\theta}$ for some $\theta >0$. According to studies by Strauss \cite{S}, Barab \cite{Ba}, Tsutsumi-Yajima \cite{TY}, and so on, there exists $ u_{\pm} \in L^2({\bf R}^n) $ such that 
\begin{align*} %\label{3}
\lim_{t \to \pm \infty} \left\| u (t,\cdot) - e^{it\Delta /2} u_{\pm} \right\|_2 = 0
\end{align*} 
holds for $2/n < \theta $ and fails for $ 0 < \theta \leq 2/n$. Therefore, in this sense, when $\theta > 2/n$ we say the nonlinearity is {\em short-range}, when $\theta \leq  {2/n}$ we say the nonlinearity is {\em long-range}, and we say the power $\theta = 2/n$ a threshold. 

Next, we consider the case where $\sigma (t) = \sigma _1 t^{-2}$ with $\sigma _1 \in (0,1/4)$. In Kawamoto-Yoneyama \cite{KY}, it was shown that the decay rate of the $L^{\infty} $--$ L^1$ dispersive estimate for free solution is of order $\CAL{O}(t^{- \theta_1 })$, $\theta _1 = n(1 + \sqrt{1-4 \sigma _1} )/4$, then the power-type nonlinearity $| u |^{\theta _2} u $ is included in the short-range class for $ \theta _1 \theta _2 > 1 $, whereas it is included in the long-range class for $\theta _1 \theta _2 = 1$ (see \cite{KM}). The key for this characterization of the threshold is the convergence and divergence of the following improper integral
\begin{align}\label{md1}
\int_1^{\infty}   \left\| F \left( u(s, \cdot ) \right)  \right\|_{\infty} ds= \int_1^{\infty} \left\| u(s, \cdot ) \right\|_{\infty}^{\theta _2} ds. 
\end{align}
Imitating the approach by Hayashi-Naumkin \cite{HN}, we can show that $\| u(s, \cdot) \|_{\infty} \leq C |s|^{-\theta _1} $ under the suitable conditions. Then we notice that if $\theta_1 \theta _2 >1$, \eqref{md1} converges and otherwise it may diverge. Consequently, imitating the approach by \cite{HN}, we can show that if \eqref{md1} converges, then the solution $u(t,x)$ converges to free solution (i.e., $\theta _2 > 1/\theta_1$ is in short-range class) and if $\theta _2 =1/\theta _1$, then $u(t,x)$ never converges to free solution (i.e., $\theta _2 = 1/ \theta _1$ is a threshold).

Now back to our model. Noting \eqref{3} and following the previous method, the critical nonlinearity can be determined through the divergence and convergence of  
\begin{align}\label{md2}
\int_1^{\infty} \left| F\left(|s|^{-n/4}(\log |s| )^{-n/2} \right)\right| ds. 
\end{align}
If we put $F(u) = |u|^{\theta_3}$. Then \eqref{md2} is equivalent to 
\begin{align}\label{md3}
\int_1^{\infty} \frac{1}{s^{n \theta _3 /4} ( \log s ) ^{n \theta _3/2}} ds. 
\end{align}
By virtue of log-decay term in \eqref{3}, improper integral \eqref{md3} is clearly integrable for $\theta _3 = 4/n$, and hence it can be expected that $4/n$ is in the short-range. On the other hand, for any $\theta _3 < 4/n$, improper integral \eqref{md3} diverges, and hence to determine the threshold for this model is difficult with only using polynomial type nonlinearities. To obtain a threshold to this model, we shall give log-modification to nonlinearity $F(u)$. Putting $F(u) = |u|^{4/n} (\log (1 + 1/|u|))^{\theta _4} $, \eqref{md2} is estimated as 
\begin{align*}
\int_1^{\infty} \left( \log \left( 1 + {s^{n/4}}{(\log s)^{n/2}} \right) \right) ^{\theta _4} \frac{ds}{s(\log s)^2} \sim 
\int_1^{\infty} \frac{ds}{s(\log s)^{2- \theta _4}},
\end{align*} 
and thus the following characterization of nonlinearities is more suitable for this model. 
\begin{Ass}\label{A2}
Let $F (u) = F_L (u)  + F_S(u) $, where 
\begin{align*}
F_S(u) = \mu_S | u |^{4/n}\left( \log \left( R + \frac{1}{|u|} \right) \right) ^{\theta}, \quad F_L(u) = \mu_L | u |^{4/n}\left( \log \left( R + \frac{1}{|u|} \right) \right) 
\end{align*}
for some $0\leq \theta <1$, $\mu_S \in {\bf R}$, and $\mu_{L} \in {\bf R}$. $F_L$ and $F_S$ are called the long- and  short-range term, respectively. Here, $R >0$ is a given constant so that for all $0 \leq \tilde{\theta} \leq 1$, and for some sufficiently small constant $0< \delta _0 \ll 1$, we have 
\begin{align*}
\inf_{t > 0} \left( 
\delta _0  \left( \log \left( R + \frac{1}{t} \right) \right)^{\tilde{\theta}}
 - \frac{\tilde{\theta} }{R t +1} \left( \log \left( R + \frac{1}{t} \right) \right)^{\tilde{\theta} -1}
\right) \geq 0 .
\end{align*}
\end{Ass}  
\begin{Rem}
Under this assumption on $R$, $t^{\delta _0} \log ( R + {1}/{t} ) ) ^{\theta} $ is a monotone increasing functions with respect to $t \geq 0$. This term appears several times, and $\delta _0 >0$ should be sufficiently small.
\end{Rem}
To state the main theorem, we define the following function spaces for $\gamma \in {\bf R}$: 
\begin{align*}
H^{\gamma ,0 } := \left\{ 
\phi \in \SCR{S}'({\bf R}^n) \, \middle| \, \| \phi  \|_{\gamma , 0} := \left( 
\int_{{\bf R}^n } \left( 1 +  |\xi|^{2} \right) ^{\gamma} \left| \SCR{F}[\phi] (\xi) \right|^2 d\xi
\right) ^{1/2} < \infty
\right\} 
\end{align*}
and 
\begin{align*}
H^{0, \gamma } := \left\{ 
\phi \in \SCR{S}'({\bf R}^n) \, \middle| \, \| \phi  \|_{ 0, \gamma } := \left( 
\int_{{\bf R}^n } \left( 1 + |x|^{2 } \right) ^{\gamma } \left| \phi (x) \right|^2 dx
\right) ^{1/2} < \infty
\right\} ,
\end{align*}
where $\SCR{F} $ denotes the Fourier transform. 
To handle log-like nonlinear terms and consider their asymptotics, the following fractional Leibniz rule is quite important: For some $\gamma >0$, we have
\begin{align} \label{8}
\left\| F_L (u) u \right\|_{\gamma , 0} \leq C \| F_L (u)  \|_{\infty} \| u \|_{\gamma , 0}. 
\end{align}
If $\gamma \in {\bf N}$, this inequality can be easily proved. However, to obtain sharp asymptotic estimates, $\gamma$ should be near $1 + 4/n$, and therefore the fractional derivative should be calculated. Roy \cite{T2, T} considered loglog nonlinear terms and established the fractional Leibniz rule for such nonlinearities. However, in \cite{T2, T}, various technical assumptions (e.g., smoothness for the nonlinearity) are made, and hence this result is difficult to apply to the present model. To overcome this difficulty, we employ the commutator estimate by Li \cite[Proposition 3.10]{L} in the following proposition, the proof of which is given in \S{2}:  
\begin{Prop}\label{P1}
Let $0 < \gamma < 1+4/n  $ for $n=1,2$ and $0< \gamma \leq 2$ for $n=3$. Then, \eqref{8} holds for smooth $u \in H^{\gamma, 0} \cap L^{\infty} ({\bf R}^n)$, and the same estimate is true if $F_L$ is replaced by $F_S$.
\end{Prop}
\begin{Rem}
Unfortunately, lack of smoothness of $F_L(u)$ for $n=3$, it is difficult to consider the case where $2 < \gamma < 1 + 4/3$. 
\end{Rem}

The aim of this paper is to investigate the asymptotic behavior for solution to \eqref{4} in $t$. To do this, we first need to show the time-in-local well-posedness for \eqref{4}, and which can be proven by imitating the proof of Theorem 3.1 in \cite{KM} (see Theorem \ref{T3}). In this process, we need to use Proposition \ref{P1} and \ref{P2}. Then we have that for all $t \in [0,T]$, there is a $C_T >0$ such that 
\begin{align*}
\left\| u(t, \cdot) \right\|_{\infty} \leq  C_T \ep ' (1 + |t |^{1/2} \log (1 + |t| ) )^{-n/2}
\end{align*}
holds. Since $C_T$ growth exponentially in $T$, this estimate does not indicate the decay estimate for $\| u(t, \cdot) \|_{\infty}$ in $t$. Hence, in order to obtain the decay estimate for $t \geq T$, we employ the approach by \cite{HN}. 

We now fix $T>0$, denote $C_T^{4/n} $ by $ \tilde{C}_1$, suppose $\ep'>0$ is sufficiently small so that $0< (\ep')^{4/3} \tilde{C}_1 (1+|\mu_L|) \ll 1$ and consider the case where $t \geq T$. By virtue of Proposition \ref{P1} and the approach by  \cite{HN}, we obtain the following theorems:
\begin{Thm}\label{T1}
Let $u_0 \in H^{\gamma ,0} \cap H^{0, \gamma}$ with $n/2 < \gamma <  1+4/n $ for $n=1,2$ and $n/2< \gamma \leq 2$ for $n=3$, and $\|  u_0\|_{\gamma , 0 } + \|  u_ 0 \|_{0, \gamma } \leq \ep ' < \ep  $, where $\ep > 0$ is sufficiently small. Then, under Assumptions \ref{A1} and \ref{A2}, there exists a unique global solution to \eqref{4} so that $u(t,x) \in C({\bf R} \, ; \, H^{\gamma ,0} \cap H^{0, \gamma })$ and 
\begin{align*}
\left\| u(t, \cdot) \right\|_{\infty} \leq C \ep ' (1 + |\zeta _2 (t)|)^{-n/2} \leq C \ep ' (1 + |t |^{1/2} \log (1 + |t| ) )^{-n/2}.
\end{align*}
\end{Thm}

\begin{Thm}\label{T2}
Under the same assumptions as in Theorem \ref{T1}, and for $u(t,x)$ as in Theorem \ref{T1}, there exist $W \in L^{\infty} ({\bf R}^n) \cap L^2 ({\bf R}^n) $ and $\tilde{C}_1 |\mu_L| (\ep ')^{4/n}  < \alpha <  \min (\gamma /2-n/4, 1)$ such that  
\begin{align*}
& \left\| \SCR{F}\left( U_0(0,t) u(t, \cdot ) \right) \mathrm{exp} \left\{i \mu_L \int_{r_0}^t F_L (|\zeta _2 (\tau)|^{-n/2} \SCR{F}(U_0(0,\tau) u(\tau, \cdot)) ) d \tau \right\} -W \right\|_{k} \\ & \leq C |\mu_L| \ep ' (\log t)^{- \alpha + \tilde{C}_1 |\mu_L| (\ep ')^{4/n} }  + C|\mu_S| \ep ' (\log t)^{ \theta -1 + \tilde{C}_1 |\mu_L| (\ep ')^{4/n}  } 
\end{align*}
if $\ep '$ is sufficiently small so that $\theta + \tilde{C}_1 |\mu_L| (\ep ')^{4/n} < 1$ and $\tilde{C}_1 |\mu_L| (\ep ')^{4/n}<  \min (\gamma /2-n/4, 1) $, where $k$ is $2$ or $\infty$. In particular there exists a function $\Phi$ such that 
\begin{align*}
& \left\| 
u(t) - e^{i F_L(|W|) (\log(\log t)) + i \Phi} U(t,0) \SCR{F}^{-1} W
\right\|_2 \\  & \leq C |\mu_L| \ep ' (\log t)^{- \alpha + \tilde{C}_1 |\mu_L| (\ep ')^{4/n} } + C|\mu_S| \ep ' (\log t)^{ \theta -1 + \tilde{C}_1 |\mu_L| (\ep ')^{4/n} } 
\end{align*}
holds.
\end{Thm}
The asymptotic behavior of long-range NLS was investigated in several studies (e.g., Strauss \cite{S}, Barab \cite{Ba}, Tsutsumi-Yajima \cite{TY}, Hayashi-Ozawa \cite{HO}, Ozawa \cite{O}, Ginibre-Ozawa-Velo \cite{GOV}, Ginibre-Velo \cite{GV}, and Hayashi-Naumkin \cite{HN}), and more rigorous analysis has been carried out in, for example, Hayashi-Naumkin-Wang \cite{HNW}, Masaki-Miyazaki-Uriya \cite{MMU}, and Okawamoto-Uriya \cite{OU}. Regarding NLS with time-dependent harmonic potentials, \cite{Ca} and Carles-Silva \cite{CS} considered general time-dependent harmonic oscillators, which include the present model, and proved various results, such as the unique existence of global-in-time solutions. Subsequently, in \cite{KM}, which focused on the case $ \zeta _2 (t) \to \infty$, it was proved that the asymptotic behavior of a solution $v(t,x)$ is 
\begin{align} \label{ad7}
\| v(t, \cdot ) \|_{\infty} \leq C (1 + |\zeta _2 (t)| )^{-n/2}
\end{align}
under suitable assumptions on $\sigma (t)$ (Assumption 1.4 in \cite{KM}). 
In particular, if $\sigma (t) t^2 \to \sigma _0 \in [0, 1/4)$, \eqref{ad7} becomes  
\begin{align*}
\| v(t, \cdot ) \|_{\infty} \leq C \ep ' |t|^{-n(1+ \sqrt{1-4 \sigma _0}) /4}
\end{align*}
(see, for instance, \cite{KY} and Geluk-Mari\'{c}-Tomi\'{c} \cite{GMT}). We notice that the threshold power of the nonlinearity is $4/(n(1+ \sqrt{1-4 \sigma _0}))$. However, the case $\sigma _0 = 1/4$, which corresponds to Assumption \ref{A1}, was not considered in \cite{KM} because of technical reasons. Accordingly, we consider this case in the present study. This result is new, and the appearance of a log nonlinearity is quite interesting (in the linear case, a similar characterization has been obtained by Ishida-Kawamoto \cite{IK}).% The log-like nonlinearities have been considered in some papers (e.g., Tao \cite{Ta}, Tristan \cite{T2,T} among others) and to consider such issues  

Main theorems can be shown by imitating the approach by \cite{HN}, herein the two norm estimates need to be shown; Lemma \ref{L4} and 
\begin{align}\label{md4}
\left\| 
U_0(0,t) u(t)
\right\|_{0, \gamma} C \ep' (\log t)^{C_{\ep'}}, 
\end{align}
where $t \geq T \gg 1$, $C>0$ and $C_{\ep'} >0$ is sufficiently small. As for the inequality \eqref{md4}, it can be shown with using \eqref{8} and the Gronwall inequality, the proof of which is stated below Lemma \ref{L4}. Proof of \eqref{8} is stated in \S{2}, herein we have to deal with the term $\| |\nabla| u \|_{\infty} $. Because $1+n/2 > 1+ 4/n $ for $n=3$, and $u$ may not be included in $H^{\rho, 0}$ with $\rho \geq 1+ 4/n$, the Sobolev embedding to $\| |\nabla| u \|_{\infty}$ fails for the case $n=3$, which demands the restriction to $\gamma$ in $n=3$. Lemma \ref{L4} can be show by \eqref{md4}, MDFM-decomposition and the same approach by \cite{HN}, the proof of which is stated in the last part of \S{3}. 
%In \S{4} we give simple examples of $\sigma (t)$ which satisfies Assumption \ref{A1}.

\section{Preliminaries}
Before proving the main theorem, we prove Proposition \ref{P1}, which is important in \cite{HN}. For $n=3$, by considering the case $\gamma =0$, $\gamma =1$ and $\gamma =2$ and interpolating them one finds Proposition \ref{P1} for $n=3$. We set $ [1+4/n]< \gamma < 1 +4/n$, and let $  \tau = \gamma  -2$ for $n=2$ and $\tau = \gamma  -4$ for $n=1$. We define 
$\tilde{F}(|u|) = |u|^{4/n} (\log (R + 1/|u|))^{\tilde{\theta}} $, $0 \leq \tilde{\theta} \leq 1$. We only consider the cases $n=2$, as the case $n=1$ is quite similar. 
\begin{Prop}\label{P2}
Let $ u \in \SCR{S}({\bf R}^n) $ with $ \SCR{F}[u] \in C_0^{\infty} ({\bf R}^n)$. Then, for $0 < \delta _2 = 1+4/n - \gamma $, we have
\begin{align} \label{11}
\left\| \tilde{F}(|u|) u \right\|_{\gamma ,0} \leq C \left\| |u|^{\delta _2} \left(\log \left(R + \frac{1}{|u|} \right) \right)^{\tilde{\theta}}\right\|_{\infty} \left\| u \right\|_{\infty}^{1 + \tau} \left\| u \right\|_{\gamma ,0}.
\end{align}
\end{Prop}
\Proof{
 We calculate $ \J{\nabla}^{\gamma } \tilde{F}(|u|) u$, where $\J{\cdot} = (1+ |\cdot|)^{1/2}$. We decompose this as $\sum \J{D_j}^{\tau} (1 +D _j ^2 ) \tilde{F}(|u|) u $, and we first calculate $D _j^2 \tilde{F} (|u|) u$, where $D_j := (- \partial _j ^2)^{1/2}$. In this process, the terms 
\begin{align}\label{16}
& \J{D_j}^{\tau} \left( \left( \log \left( R + \frac{1}{|u|}\right) \right)^{\tilde{\theta} -1} + \left( \log \left( R + \frac{1}{|u|}\right) \right)^{\tilde{\theta}-2} \right)  \left( (u_j)^2 + u_{jj}u \right) \CAL{O}(|u|^{4/n-1}) \\ & \quad + \mbox{ (similar terms )} \nn
\end{align} 
and 
\begin{align}\label{17}
& \J{D_j}^{\tau} \left( \log \left( R + \frac{1}{|u|}\right) \right)^{\tilde{\theta}} (u_{jj} u + (u_j) ^2) \CAL{O}(|u|^{4/n -1})   + \left( \mbox{ similar terms } \right) %\nn
\end{align}
appear. In the followings, we only estimate the terms in \eqref{17}, because the terms in \eqref{16} can be handled by the same approach for to handle \eqref{17}, and find their $L^2$ norm can be estimated as $ C \| u \|_{\infty}^{4/n} \| u\|_{\gamma  ,0}$ because $(\log (R + 1/t))^{\tilde{\theta} -1} $ and $(\log (R + 1/t))^{\tilde{\theta} -2} $ are uniformly bounded for all $t \geq 0$.

The term associated with $u_{jj} u$ is easy to estimate compared with the term $u_j^2$, and hence we only consider the latter. By using the identity 
\begin{align*}
\frac{\bar{u}}{|u|} ( u_j)^2 |u|^{- \beta } + |u|^{1- \beta} u_{jj} + \frac{\beta}{1-\beta} |u|^{1- \beta} u_{jj} = \frac{\bar{u}}{|u|} \partial _j \left( 
\frac{u}{|u|^{\beta}} u_j
\right) + \frac{\beta}{1- \beta} \partial _j (|u|^{1- \beta} u_j) 
\end{align*}
with $\beta = 1- \tau$, we obtain 
\begin{align*}
u_j^2 \sim |u|u_{jj} + |u|^{1- \tau } \partial _j (\CAL{O}(|u|^{\tau}) u_j), 
\end{align*}
where $\SCR{F}[u] \in C_0^{\infty} ({\bf R}^n)$ ensures the boundedness of $|u|^{1-\tau} \partial _j (\CAL{O}(|u|^{\tau}) u_j)$. Therefore, we estimate 
\begin{align*}
 D_j^{\tau} \left( \log \left( R + \frac{1}{|u|}\right) \right)^{\tilde{\theta}} \CAL{O}(|u|^{4/n -\tau }) \partial _j (\CAL{O}(|u|^{\tau}) u_j) . 
\end{align*}
Let $f =  ( \log ( R + {1}/{|u|}) )^{\tilde{\theta}} \CAL{O}(|u|^{4/n -\tau })  $ and $g = D_j^{\tau} (\CAL{O}(|u|^{\tau}) u_j) $, i.e, $\partial _j (\CAL{O}(|u|^{\tau}) u_j ) = \partial_j D_j^{-1} \cdot D_j^{1- \tau} g $. Then, by the sharp  commutator estimate by Li \cite[Proposition 3.10]{L}, we obtain 
\begin{align*}
\left\| 
D_j^{\tau} (f D_j^{ 1- \tau} (\partial _j D_j^{-1} g) ) - f \partial_j g
\right\|_2 \leq C \left\| \partial_j f \right\|_{\infty} \left\| \partial _j D_j^{-1} g \right\|_2.
\end{align*}
Hence, we have  
\begin{align*}
\left\| 
D_j^{\tau} (f \partial _j (\CAL{O}(|u|^{\tau}) u_j)  ) 
\right\|_2 &\leq \left\| f \partial _j g \right\|_2 + C \left\| \partial _j  f \right\|_{\infty} \left\| g \right\|_{2} \\ & \leq  \left\| f  \partial _j g \right\|_2 + C \left\| \partial _j f \right\|_{\infty} \left\| D_j^{\tau} \CAL{O}(|u|^{\tau}) u_j \right\|_{2}. 
\end{align*}
Regarding $\left\| f \partial _j g \right\|_2 $, we first estimate it as 
\begin{align*} 
\left\| f \partial _j  g \right\|_2  \leq \left\| \CAL{O}(|u|^{4/n-\tau -1}) (\log (R + 1/|u|) )^{\tilde{\theta}}\right\|_{\infty} \left\| u \partial _j g \right\|_2. 
\end{align*}
Subsequently, we use Proposition 3.10 in \cite{L} and obtain 
\begin{align*}
 \left\| u \partial _j g \right\|_2 \leq \left\| D_j ^{\tau} uD_j^{1-\tau} ( \partial _j D_j^{-1} g)  \right\|_2 + C \left\| \partial _j u \right\| _{\infty} \left\| D_j^{\tau} \CAL{O}(|u|^{\tau} ) u_j \right\|_{2}. 
\end{align*} 
As
\begin{align*}
D_j^{\tau} u D_j^{1-\tau} ( \partial _j D_j^{-1} g)  = - D_j^{\tau} \CAL{O} (|u|^{\tau}) u_j^2  + D_j^{\tau} \partial _j \left( \CAL{O}(|u|^{ \tau}) u u_j \right),
\end{align*}
for $\delta _2 = 4/n - \tau -1$, the Gagliardo--Nirenberg interpolation inequality implies
\begin{align}\label{12}
\left\| f \partial _j g \right\|_2  \leq C \left\| |u|^{\delta _2} \left(\log \left(R + \frac{1}{|u|} \right) \right)^{\tilde{\theta}}\right\|_{\infty} \left\| u \right\|_{\infty}^{1 + \tau} \left\| u \right\|_{\gamma ,0}.
\end{align}
Moreover, as 
\begin{align*}
\left\| \partial _j f \right\|_{\infty} &\leq C \left\| \left(\log \left(R + \frac{1}{|u|} \right) \right)^{\tilde{\theta}} \CAL{O}(|u|^{\delta _2} ) \partial _j u \right\|_{\infty} + (\mbox{ similar terms })\\ & \leq C  \left\| |u|^{\delta _2} \left(\log \left(R + \frac{1}{|u|} \right) \right)^{\tilde{\theta}} \right\|_{\infty} \left\| \partial_j u \right\|_{\infty},
\end{align*}
we have 
\begin{align}\label{13}
\left\| \partial _j f \right\|_{\infty} \left\|g \right\|_{2} \leq C \left\| |u|^{\delta _2} \left(\log \left(R + \frac{1}{|u|} \right) \right)^{\tilde{\theta}}\right\|_{\infty} \left\| u \right\|_{\infty}^{1 + \tau} \left\| u \right\|_{\gamma ,0}. 
\end{align}
By combining \eqref{12} and \eqref{13}, Proposition \ref{P2} is proved.

}

As $t^{\delta _2} (\log (R + 1/t))^{\tilde{\theta}} $ is a monotone increasing function, the assumption $\| u \|_{\infty} \leq C $ implies 
\begin{align*}
\left\| |u|^{\delta _2} \left(\log \left(R + \frac{1}{|u|} \right) \right)^{\tilde{\theta}}\right\|_{\infty} \leq  \|u \|^{\delta _2}_{\infty} \left(\log \left(R + \frac{1}{\|u \|_{\infty}} \right) \right)^{\tilde{\theta}}. 
\end{align*}
Therefore, by Proposition \ref{P2}, we obtain 
\begin{align}\label{14}
\left\|  \tilde{F}(|u|) u \right\|_{\gamma ,0} \leq C \left\|  u \right\|_{\infty}^{4/n}\left(\log \left(R + \frac{1}{\|u \|_{\infty}} \right) \right)^{\tilde{\theta}} \left\| u \right\|_{\gamma ,0} = C \tilde{F}(\| u \|_{\infty}) \| u \|_{\gamma,0}. 
\end{align} 
As $\tilde{F}$ is also monotone increasing, we have 
\begin{align}\label{15}
\left\| \tilde{F} (|u|)  \right\|_{\infty}  = \tilde{F}(\| u \|_{\infty}) . 
\end{align} 
Combining \eqref{14} and \eqref{15}, one can prove Proposition \ref{P1}. 

Thanks to the proposition \ref{P1} and interpolation theorem, we find the following Proposition.
\begin{Prop}\label{P2}
Let $u,v \in \SCR{S} ({\bf R }^n)$ with $\SCR{F}[u], \SCR{F}[v] \in C_0^{\infty} ({\bf R}^n)$. Then, for the same $\gamma$ in Theorem \ref{T1}, 
\begin{align}\label{18}
\left\| 
\tilde{F}(|u|) u - \tilde{F}(|v|) v
\right\|_{\gamma ,0} \leq C \left( \left\| \tilde{F}(|u|) \right\|_{\infty} + \left\| \tilde{F} (|v|) \right\|_{\infty} \right) \left\| u-v \right\|_{\gamma ,0}.
\end{align}
\end{Prop}
~~ \\ 
By Proposition \ref{P1} and \ref{P2}, one can prove the existence of a unique local-in-time solution to \eqref{4}; 
\begin{Thm} \label{T3}
Let $T \gg r_0$ be fixed and $u_0 \in H^{\gamma , 0} \cap H^{0, \gamma} $ with the same $\gamma$ in Theorem \ref{T1}. We assume that $\| u_0 \|_{\gamma ,0} + \|  u_0 \|_{0, \gamma } \leq \ep ' \leq \ep = \ep (T) $, where $\ep (T) >0$ is sufficiently small compared with $T$. Then, under Assumptions \ref{A1} and \ref{A2}, there exist a unique local-in-time solution to \eqref{4} and $C_{T} >0$ so that $u(t,x) \in C( [-T, T] \, ; \, H^{\gamma ,0} \cap H^{0, \gamma})$, $C_{T} \ep (T) < 1 \to 0$ as $\ep (T) \to 0$, 
\begin{align*}
\sup_{t \in [-T, T]} \left\| 
u(t, \cdot)
\right\|_{\infty} \leq C_{T} \ep ' (1 + |\zeta _2 (t)|)^{-n/2},
\end{align*} 
and 
\begin{align*}
\sup_{t \in [-T, T]}\left\| 
U_0(0, t) u(t, \cdot)
\right\| _{\gamma ,0 } \leq C_{T} \ep '.
\end{align*}
 
\end{Thm}
The proof of this theorem is similar to that of Theorem 3.1 in \cite{KM}, and thus we omit it. 

By this theorem, we have $\left\| 
u(s, \cdot)
\right\|_{\infty} \leq C_{r_0} \ep '$ and $\left\| 
U_0(0, t) u(t, \cdot)
\right\| _{\gamma ,0 } \leq C_{r_0} \ep '$ for all $t \in [-r_0,r_0]$ and for some $C_{r_0 } >0$. As $T$ is sufficiently large compared with $r_0$, one can assume that $0< C_{r_0} \ll C_T$, and hence without loss of generality, $C_{r_0}$ can be denoted by $C$.

\section{Proof of Theorems \ref{T1} and \ref{T2}} 
Herein, we prove Theorems \ref{T1} and \ref{T2}. Theorem \ref{T2} can be proved by the arguments in the proof of Theorem \ref{T1} and the approach in \cite{HN}. Hence, we only prove Theorem \ref{T1}. We first show that $C_T$ in Theorem \ref{T3} can be selected independently of $T$; subsequently, we show Theorem \ref{T1} by a continuation argument. The second part immediately follows from the first; hence, it suffices to show the first part.  

To simplify the proof, we introduce the following lemma, which can be proved by using Lemma 2.11 in \cite{KM} and the arguments in the proof of Theorem 4.2 in \cite{KM}.
\begin{Lem}\label{L4}
For some smooth $u(t)$, we have the estimates 
\begin{align*}
\left\| 
u(t)
\right\|_{\infty} \leq C |\zeta _2 (t)|^{-n/2} \| \SCR{F} U_0(0,t) u(t) \|_{\infty} + C |\zeta _2 (t)|^{-n/2} \left| 
\frac{\zeta _1 (t)}{ \zeta _2 (t)}
\right|^{\alpha} \| U_0(0,t) u(t) \|_{0, \gamma}  
\end{align*}
and 
\begin{align} \label{ad1}
\left\| 
\SCR{F} U_0(0,t) u(t)
\right\|_{\infty} \leq C \ep ' + C \int_{r_0}^t \left( |\mu_L|  I_1 (\tau)  + |\mu_S |  I_2 (\tau)  \right) d \tau
\end{align}
with 
\begin{align*}
I_1 (t) &= \left| 
\frac{\zeta _1 (t)}{ \zeta _2 (t)}
\right|^{\alpha} F_L \left( |\zeta _2 (t)|^{-n/2} \left\| U_0(0,t)u(t) \right\|_{0, \gamma '} \right)  \left\| U_0(0,t)u(t) \right\|_{0, \gamma }
\end{align*}
and 
\begin{align*}
I_2 (t) = F_S \left( |\zeta _2 (t)|^{-n/2} \left\| U_0(0,t)u(t) \right\|_{0, \gamma '} \right)  \left\| U_0(0,t)u(t) \right\|_{0, \gamma ' },
\end{align*}
where $t > r_0$, $\gamma '> n/2$, $\gamma > n/2 + 2 \alpha$, and $0< \alpha \leq 1$. 
\end{Lem}
~~ \\ 

We now prove Theorem \ref{T1}. We first estimate the term $\| U_0(0,t) u(t)\|_{0, \gamma} $. As 
\begin{align*}
\left\| 
u(t)
\right\|_{\infty} \leq C_T \ep ' \left( 1 + |\zeta _2 (t)| \right)^{-n/2} 
\end{align*}
for $t \in [-T,T]$, and $F_L (t)$ and $F_S (t)$ are monotone increasing functions on $t \geq 0$, 
we have 
\begin{align*}
& \left\| 
U_0(0,t) u(t)
\right\|_{0, \gamma} \\ & \leq \| u_0 \|_{0, \gamma} + C \int_{0}^t \left( 
|\mu_L| F_L (\| u(s) \|_{\infty} ) + 
|\mu_S| F_S (\| u(s) \|_{\infty} ) 
\right)  \left\| U_0(0,s) u(s) \right\|_{0,\gamma} ds .
\end{align*}
Then, by the inequality 
\begin{align*}
C \int_{0}^{r_0} \left( 
|\mu_L| F_L (\| u(s) \|_{\infty} ) + 
|\mu_S| F_S (\| u(s) \|_{\infty} ) 
\right)  \left\| U_0(0,s) u(s) \right\|_{0,\gamma} ds \leq C \ep ' 
\end{align*}
we have 
\begin{align*}
& \left\| 
U_0(0,t) u(t)
\right\|_{0, \gamma} \\ & \leq 
C\ep ' + C C_T^{4/n} \int_{r_0}^t (\ep ')^{4/n} (1 + |\zeta _2 (s)| )^{-2} \Bigg( 
| \mu_L  | \log \left( R + \frac{(1 + | \zeta _2 (s) | )^{n/2}}{\ep' C_T } \right) \\ & \quad + 
| \mu_S  | \left( \log \left( R + \frac{(1 + |\zeta _2 (s) | )^{n/2}}{\ep ' C_T } \right) \right)^{\theta}
\Bigg)  \left\| U_0(0,s) u(s) \right\|_{0,\gamma} ds \\ 
 & \leq 
C \ep ' + C C_T^{4/n}(\ep ')^{4/n-1} \int_{r_0}^t  (1 + |\zeta _2 (s)|)^{-2} \Bigg( 
| \mu_L |  \ep ' \left(  \log ( R \ep' C_T + 1 + | \zeta _2 (s)|^{n/2} ) + \log (\ep ' C_T) \right) \\ & \quad + 
| \mu_S | \ep '  \left(  \left( \log ( R \ep ' C_T + 1 + | \zeta _2 (s) |^{n/2}) \right)^{\theta} + \left( \log (\ep ' C_T) \right)^{\theta}\right)
\Bigg)  \left\| U_0(0,s) u(s) \right\|_{0,\gamma} ds. 
\end{align*}
Clearly, for $0 < \tilde{\theta} \leq 1$, we have 
\begin{align*}
& \ep' \left(  \left( \log ( R \ep' C_T + 1 + |\zeta _2 (s)|^{n/2})  \right)^{\tilde{\theta}} + \left(  \log (\ep ' C_T)  \right)^{\tilde{\theta}} \right) \\ & \leq C \left( C_T^{-1} \ep' C_T \left( \log (\ep' C_T) \right)^{\tilde{\theta}} + \ep' \left( \log |\zeta _2 (s)| \right)^{\tilde{\theta}}  \right) \\ & \leq C \left( C_T^{-1}  + \ep' ( \log s )^{\tilde{\theta}} \right), 
\end{align*}
where we use $|\log (s^{1/2} \log s)| \leq C(\log s + \log \log s) \leq C \log s$ for $s > r_0 \gg 1$. Then, the inequalities $(1 + |\zeta _2 (s) | )^{-2} \leq C s^{-1} (\log s)^{-2}$ and ${\displaystyle \int_{r_0}^t (\cdots) ds \leq C \ep ' + \int_{0}^t (\cdots) ds} $ imply 
\begin{align*}
& \left\| 
U_0(0,t) u(t)
\right\|_{0, \gamma} \\ & \leq  C\ep'+ C C_T^{4/n-1} (\ep')^{4/n -1} (| \mu_S | + | \mu_L | ) \int_{0}^t (1+s)^{-1} (\log (2+s))^{-2} \left\| U_0(0,s) u(s) \right\|_{0,\gamma}  ds  \\ & \quad + 
C C_T^{4/n}(\ep')^{4/n} \int_{0}^t (1+s)^{-1} \left( 
| \mu_L |  (\log (2+s) )^{-1} + | \mu_S | (\log (2+s))^{-2 + \theta}
 \right)\left\| U_0(0,s) u(s) \right\|_{0,\gamma}  ds. 
\end{align*}
Hence, 
\begin{align} \nn 
\left\| 
U_0(0,t) u(t)
\right\|_{0, \gamma} & \leq C \ep '  \mathrm{exp} \left\{  C_T^{4/n-1} (\ep')^{4/n-1} ( |\mu_S| + |\mu_L| )  +  C_T^{4/n}(\ep')^{4/n} (|\mu_S| + |\mu_L| \log (\log t)   \right\} \\ & \leq 
C \ep ' \mathrm{exp} \left\{ C_{0,\ep '} \left( | \mu _S | +   |\mu_L | \log (\log t) \right) \right\},  \label{ad2}
\end{align}
where $C_{0, \ep '} = C_T^{4/n} (\ep')^{4/n } $ and  $0<C_{0, \ep '} \ll 1$.
\\ ~~ \\ 
{\bf Estimation of $\| u(t) \|_{\infty} $.} \\ 
By Lemma \ref{L4}, we have 
\begin{align*}
\left\| u(t) \right\|_{\infty} & \leq C \ep ' | \zeta _2 (t) |^{-n/2} (\log t)^{ - \alpha} 
 \mathrm{exp} \left\{  C_{0, \ep '} (| \mu_S | + | \mu_L | \log (\log t)  ) \right\} \\ & \quad 
 + C | \zeta _2 (t) |^{-n/2} \left\| \SCR{F} U_0(0,t) u(t) \right\|_{\infty}. 
\end{align*}
Let now 
\begin{align*}
v (t) := U_0(0,t) u(t) 
\end{align*}
and $\gamma > \gamma '$. Then, $\| \cdot \|_{\gamma '} \leq \| \cdot  \|_{\gamma}$ and \eqref{ad1} imply 
\begin{align*}
\left\| 
\SCR{F} U_0(0,t) u(t)
\right\| _{\infty}  & \leq C\ep ' 
+ \int_{r_0}^t \Big\{ | \mu_{L} |  (\log s)^{ - \alpha} F_L (\zeta _2 (s)^{-n/2}  \| v(s) \|_{0, \gamma} ) \| v(s) \|_{0, \gamma} \\ & \quad + |\mu_{S} | F_S (\zeta _2 (s) ^{-n/2} \| v(s) \|_{0, \gamma} ) \| v(s) \|_{0, \gamma} \Big\} ds . 
\end{align*}
If $f(t) = t^{4/n} (\log (R + 1/t))^{\tilde{\theta}} $ with $0 \leq \tilde{\theta} \leq 1$ (a monotone increasing function for $t \geq 0$) and $C_{1} (\ep ', \mu_L) := 4 (C_{0, \ep '} ) |\mu_L|/n$, then \eqref{ad2} implies 
\begin{align*}
& f(|\zeta _2 (s)|^{-n/2} \| v(s) \|_{0, \gamma} ) \\  & \leq C s^{-1} (\log s)^{-2} (\ep ')^{4/n} (\log s)^{C_1(\ep ', \mu_L)} \left\{ 
 \left( 
\log \left( 
R + \frac{s^{n/4} (\log s)^{n/2} }{C \ep' (\log s)^{ C_{0, \ep '} |\mu_L |} }
\right) 
\right) ^{\tilde{\theta}}
\right\} \\ & \leq 
C  s^{-1} (\log s)^{-2 + C_1 (\ep', \mu_L)} \times
\begin{cases}
(\ep') ^{4/n -1}, & (s (\log s)^{2-C_1(\ep',\mu_L)})^{n/4} \leq 1/( 2C\ep ') , \\ 
( \ep ') ^{4/n} (\log s)^{\tilde{\theta}},  & (s (\log s)^{2-C_1(\ep',\mu_L)})^{n/4} \geq 1/( 2C\ep ') ,
\end{cases} 
\\ & \leq 
C s^{-1} (\log s)^{-2 + C_1 (\ep', \mu_L) + \tilde{\theta}} \ep ^{4/n -1}. 
\end{align*} 
Therefore,
\begin{align*}
\left\| 
\SCR{F} U_0(0,t) u(t)
\right\| _{\infty}  & \leq C \ep '   + C |\mu_L| (\ep ')^{4/n} \int_{r_0}^t s^{-1} (\log s)^{-1 - \alpha + C_1 (\ep' , \mu_L) + C_{0, \ep '} | \mu_L | } \\ & \qquad \quad 
+ C |\mu_S| (\ep ')^{4/n} \int_{r_0}^t s^{-1} (\log s)^{- 2 + \theta + C_1 (\ep', \mu_L) + C_{0, \ep '} | \mu_L | }  \\ & \leq C \ep ',
\end{align*}
where $\tilde{C}_{1} ( \ep ', \mu_L)  :=  C_1 (\ep ' , \mu_L) + C_{0, \ep '} | \mu_L |$, and we assume that $- \alpha +\tilde{C}_1 (\ep, \mu_L) < 0 $ and $\theta + \tilde{C}_1 (\ep, \mu_L)   < 1$, which holds if $\ep '$ is sufficiently small because $\tilde{C}_{1} ( \ep ', \mu_L) \to 0  $ as $\ep ' \to 0$.    

\section{Models of $\sigma (t)$ satisfying Assumption \ref{A1}} 
Herein, we provide simple examples of $\sigma (t)$ satisfying Assumption \ref{A1}. For $\alpha \neq 0$, let 
\begin{align*}
\sigma (t) = 
\begin{cases}
\alpha ^2, & |t| \leq r_0, \\ 
t^{-2 }/4, & |t| > r_0. 
\end{cases}
\end{align*}
Then, $\zeta _j (t)$ in \eqref{1} can be written as 
\begin{align*}
\zeta _1 (t) = \cos (\alpha t), \quad \zeta_2 (t) = \alpha ^{-1} \sin (\alpha t)
\end{align*}
for $|t| \leq r_0$ and 
\begin{align*}
\zeta _1 (t) = c_{11, \pm} y_1 (t) + c_{12, \pm} y_2 (t), \quad 
\zeta _2 (t) = c_{21, \pm} y_1 (t) + c_{22, \pm} y_2 (t), 
\end{align*}
where $c_{ij, \pm} \in {\bf R}$ are some constants, $y_1 (t) = |t| ^{1/2}$, and $y_2 (t) = |t|^{1/2} \log |t|$. By the continuation condition for $\zeta _j (t)$ and $\zeta _j '(t)$, we notice 
\begin{align} \nn 
c_{11, \pm} r_0 ^{1/2} + c_{12, \pm} r_0 ^{1/2} \log r_0 &= \cos (\alpha (\pm r_0) ), \\ 
c_{21, \pm} r_0 ^{1/2} + c_{22, \pm} r_0 ^{1/2} \log r_0 &= \alpha ^{-1}  \sin (\alpha (\pm r_0) )
\label{ad4}
\end{align}
and 
\begin{align}\nn
c_{11, \pm} + c_{12, \pm} \log r_0 + 2 c_{12, \pm} &= \mp  2 \alpha r_0 ^{1/2} \sin (\alpha (\pm r_0)), \\ 
c_{21, \pm} + c_{22, \pm} \log r_0 + 2 c_{22, \pm} & = \pm 2r_0^{1/2} \cos (\alpha (\pm  r_0 )).
\label{ad5}
\end{align}
Condition \eqref{ad3} implies $c_{12, \pm} = 0$ and $c_{22, \pm} \neq 0$; therefore, 
\begin{align*}
c_{11, \pm} = r_0^{-1/2} \cos (\alpha r_0) = - 2 \alpha r_0 ^{1/2} \sin (\alpha r_0), 
\end{align*}
which implies 
\begin{align} \label{ad8}
\alpha r_0 \tan (\alpha r_ 0 ) = - 2.
\end{align}
Moreover, by \eqref{ad4} and \eqref{ad5}, we obtain 
\begin{align} \label{ad9}
2c_{22, \pm} = \pm 2 r_0^{1/2} \cos (\alpha r_0) \mp  \alpha ^{-1} r_0^{-1/2} \sin (\alpha r_0) = \frac{\pm 2 \cos (\alpha r_0)}{\alpha r_0 ^{1/2}} \left( 
\alpha r_0  - \tan (\alpha r_0)
\right) \neq 0.
\end{align}
Hence, if $\alpha r_0 $ satisfies \eqref{ad8} and \eqref{ad9}, then $\sigma (t) $ in \eqref{1} satisfies Assumption \ref{A1}.
\\ ~~ \\ 
{\bf \large Acknowledgement} \\ 
The author is supported by the Grant-in-Aid for Young Scientists \#20K14328 from JSPS.
%\appendix
%\def\thesection{Appendix \Alph{section}}
%\section{Proof of Theorem \ref{}}

\end{document}